\documentstyle[12pt,amssymbt]{amsart}
\renewcommand{\a}{\alpha}
\renewcommand{\b}{\beta}

\renewcommand{\d}{\delta}
\newcommand{\e}{\varepsilon}
\newcommand{\z}{\zeta}
\newcommand{\vt}{\vartheta}
\renewcommand{\l}{\lambda}
\newcommand{\s}{\sigma}
\renewcommand{\t}{\tau}

\newcommand{\f}{\varphi}

\renewcommand{\L}{\Lambda}

\renewcommand{\O}{\Omega}

\newcommand{\A}{{\cal A}}

\newcommand{\p}{{\cal P}}

\newcommand{\C}{\Bbb C}
\newcommand{\T}{\Bbb T}
\newcommand{\pp}{{\Bbb P}}
\newcommand{\dd}{{\Bbb D}}

\newcommand{\0}{{\Bbb O}}
\newcommand{\bO}{\{ {\Bbb O} \}}

\newcommand{\df}{\stackrel{\rm{def}}{=}}

\newcommand{\Ker}{\operatorname{Ker}}

\newcommand{\beq}{\begin{equation}}
\newcommand{\eeq}{\end{equation}}
\newcommand{\bay}{\begin{eqnarray}}
\newcommand{\ey}{\end{eqnarray}}
\newcommand{\ra}{\rightarrow}

\newcommand{\be}{\infty}

\newcommand{\bl}{\blacksquare}

\renewcommand{\pf}{{\bf Proof. }}

\newtheorem{thm}{\hspace{\parindent}Theorem}[section]

\newtheorem{cor}[thm]{\hspace{\parindent}Corollary}
\newtheorem{lem}[thm]{\hspace{\parindent}Lemma}

\date{\mbox{\today}}
\begin{document}

\overfullrule=5pt

\newcommand{\vse}{\vspace{.2in}}
\newcommand{\N}{\Bbb N}
\newcommand{\U}{\Upsilon}
\renewcommand{\theequation}{\thesection.\arabic{equation}}

\newcommand{\ov}{\overline}
\newcommand{\shd}{\hspace{-2.5ex plus  -.2ex}.\hspace{1.4ex plus .1ex}}

%\font\Bbb=msbm10 scaled \magstep1

\newcommand{\thd}{\hspace{-1.1ex}{\bf.}\hspace{1.1ex}}

\title{ Continuity Properties of Best Analytic Approximation}
\author{V.~V.~Peller}
 \address{V.~V.~Peller\\
Department of Mathematics \\ 
Kansas State University \\
Manhattan \\
Kansas 66502 }

\author{N.~J.~Young}
\address{N.~J.~Young\\
Department of Mathematics and Statistics\\
Lancaster University\\
Lancaster LA1 4YF \\
England}

\thanks{ V.V.Peller's research was supported by an NSF grant in Modern Analysis. 
N.J.Young wishes to thank
the Mathematics Departments of Kansas State University and
the University of California at San Diego, and the Mathematical Sciences Research
Institute for hospitality while this work was carried out.  Research at MSRI is
supported in part by NSF grant DMS-9022140. }

\subjclass{AMS Subject Classifications: 30E10, 47B35, 93B36.}

\setcounter{equation}{0}
\setcounter{section}{0}

\maketitle

\begin{abstract}
Let $\A$ be the operator which assigns to each  $m \times n$ matrix-valued
 function on
the unit circle with entries in  $H^\infty + C$ its
unique superoptimal approximant in the space of bounded analytic
 $m \times n$ matrix-valued functions in the open unit disc.
We study the continuity of $\A$ with respect to various norms.  Our main
result is that, for a class of norms satifying certain natural axioms, $\A$
is continuous at any function whose superoptimal singular values are
non-zero and is such that certain associated integer indices are equal to 1.
We also obtain necessary conditions for continuity of $\A$ at point and a
sufficient condition for the continuity of superoptimal singular values.

\end{abstract}

\section{Introduction}

The problem of finding a best uniform approximation of a given bounded function on the
unit circle by an analytic function in the unit disc is a natural one from the viewpoint
of pure mathematics and it also has engineering applications, for example in $H^\be$
control [F], broadband impedance matching [He] and robust identification [Par].  In
these contexts, to effect a design or construct a model, one must compute such
a best approximation, and in order that numerical computations have validity it is
important that the solution to be computed depend continuously on the input data, for
otherwise the imperfect precision of floating point arithmetic may lead to highly
inaccurate results.  It is therefore somewhat disconcerting that, with respect to the
$L^\be$ norm, the operator of best analytic approximation is discontinuous everywhere
except at points of $H^\be$ 
[M, Pa].  Nevertheless engineers regularly compute such approximations and appear to
find the results reliable.  A way to account for this would be to show that best
analytic approximation {\em is} continuous on suitable Banach subspaces of $L^\be({\Bbb
T})$ with norms which majorise the uniform norm, or at least, is continuous at most
points of the space.  One can expect that most functions of
engineering interest will lie in one of these well-behaved subspaces, and that the
errors introduced by computer arithmetic will result in perturbations which are small in
the associated norm.  We are thus led to ask for which Banach spaces $X \subset 
L^\be({\Bbb T})$ the operator $\A$ of best analytic approximation maps $X$ into $X$ and is
continuous at a generic point of $X$ (in some sense).  This question has been thoroughly
analysed for the case of scalar-valued functions.  It was shown in [P1] that, for spaces 
$X \subset H^\infty+C$ satisfying some natural axioms, the restriction of $\A$
to $X$ is continuous with respect to the norm of $X$ at a function $\f$
if and only if $\| H_\f \|$ is a simple singular value of the Hankel
operator $ H_\f$.

Analogous questions for matrix-valued functions are also of interest,  
particularly for their relevance to engineering applications.  They are a good
deal more complicated than in the scalar case.  To begin with, there is
typically no unique best analytic approximation in the matrix case, when we
measure closeness by the  $L^\be$ norm.  In order to specify an
approximation uniquely and so obtain a well formulated question of
continuity we can use a more stringent criterion of approximation.  The
notion of a {\em superoptimal} approximation is a natural one for
matrix-valued functions: by imposing the condition of the minimisation of
the suprema of all singular values of the error function it gives a unique
best approximant in many cases.  Here is a precise definition.

Denote by $M_{m,n}$ the space of $m \times n$ complex matrices endowed with the
operator norm as a space of linear operators from $\;\,\C^n$ to
$\;\,\C^m$ with their
standard inner products.  Let $H^\infty (M_{m,n})$ denote the space of bounded  
analytic $M_{m,n}$-valued functions on the unit disc ${\Bbb D}$ with supremum norm:
$$
||Q||_{H^\infty} \stackrel{{\rm def}}{=}
||Q||_\infty \stackrel{{\rm def}}{=} \sup_{z \in {\Bbb D} }
||Q(z)||_{M_{m,n}}.   
$$
Similarly, $L^\infty (M_{m,n})$ denotes the space of essentially bounded
Lebesgue measurable $M_{m,n}$-valued functions on $\T$ with essential
supremum norm.  By Fatou's theorem [H, p.34] functions in $H^\infty (M_{m,n})$
have radial limits a.e. on $\T$, so that $H^\infty (M_{m,n})$ can be
embedded isometrically in $L^\infty (M_{m,n})$, and we shall often tacitly
regard elements of $H^\infty (M_{m,n})$ as functions on the unit circle. Where  
there is no risk of confusion we shall sometimes write
$H^\infty,\,\; L^\infty$  for $ H^\infty (M_{m,n}),\;\, L^\infty (M_{m,n}).$ We define
$H^\infty + C$ to be the space of (matrix-valued) functions on
$\T$ which are expressible as the sum of an $H^\infty$ function and a
continuous function on $\T$.
For any matrix $A$ we denote the transpose of $A$ by $A^t$ and the 
singular values or $s$-numbers of $A$ by
$$
s_0 (A) \geq s_1 (A) \geq \dots \geq 0.
$$
For $F \in L^\infty (M_{m,n})$ we define, for $j = 0, 1, 2, \dots,$
$$                                                                              
s_j^\infty (F) \stackrel{{\rm def}}{=} {\rm ess}\;\sup_{|z|=1}s_j(F(z))  
$$
and
$$
s^\infty (F) \stackrel{{\rm def}}{=}
(s^\infty_0 (F), s^\infty_1 (F), s^\infty_2 (F), \dots).
$$
We shall say that $Q \in H^\infty(M_{m,n})$ is a superoptimal $H^\infty$
approximant to 
$\Phi \in L^\infty(M_{m,n})$ if $s^\infty(\Phi-Q)$ is a minimum over
$Q \in H^\infty$ with respect to the lexicographic ordering. 

It was proved in [PY1] that if an $m \times n$ matrix function $\Phi$
 is in $ H^\infty +C$ then there is a
 unique superoptimal approximant to $\Phi$ in  $H^\infty (M_{m,n})$.
 We shall denote this approximant  by ${\cal A}\Phi$.  In [PY1], in addition to
proving uniqueness, we obtained detailed structural information about the
``superoptimal error" $\Phi-{\cal A}\Phi$ and we established several heredity
results (that is, theorems of the form $``\Phi \in X$ implies $ {\cal A}\Phi \in
X$" for various function spaces $X$).  In any space which does have this
heredity property it is natural to ask whether $\A$ acts continuously. 
We shall show that for a substantial class of norms there are many
continuity points of $\A$.  We cannot, however, expect $\A$ to be continuous everywhere:
it is shown in [P1] that, for scalar functions,
 $\A$ is discontinuous with respect to virtually any norm at every
$\f$ for which $\|H_\f \|$ is a multiple singular value of $\| H_\f \|$, and it follows
that (matricial) $\A$ is discontinuous at the matrix function diag$\{\f, 0,\cdots \}$.

We shall study spaces  $X \subset L^2({\Bbb T})$ of functions for which the following
axioms hold.  Denote by $\pp_{+},~\pp_{-}$ the orthogonal projections from $L^2({\Bbb T})$
onto the Hardy space $H^2$ and its orthogonal complement $H^2_{-}$ in $L^2({\Bbb T}).$
For a space $X \subset L^2({\Bbb T})$ we denote by $X_{+}$ the space
$\{\pp_{+} f: f \in X\}$ and by $X_{-}$ the space $\{\pp_{-}f: f \in X\}$.
The axioms are:
\newline (A1)  If $f \in X$ then $\bar{f} \in X$ and $\pp_{+}f \in X$;
\newline (A2) $X$ is a Banach algebra with respect to pointwise multiplication;
\newline (A3) the set of trigonometric polynomials is dense in $X$;
\newline (A4) every multiplicative linear functional on $X$ is of the form $f \mapsto
f(\z)$ for some $\z \in \T$;
\newline (A5) if $f \in X_{+}$ and $h \in H^\be$ then $\pp_{+}(\bar{h}f) \in X_{+}$.

The following fact is well known.

\begin{lem}\thd
 $ X_{+}$ with the restriction of $\|\cdot\|_X$
 is a commutative Banach algebra
 whose maximal ideal space is the closed unit disc {\rm clos} $\Bbb D$.
\end{lem}

\pf  By the Closed Graph Theorem $\pp_{+}$ is continuous on $X$, and so its 
range $X_{+}$ is a closed subspace of $X$.  Functions in $X_{+}$ are continuous
on $\T$  (the Gelfand topology of $X$ on $\T$ is compact and refines the
 natural
 topology, hence coincides with it), and their negative Fourier coefficients 
vanish.  Hence $X_{+} \subset A({\Bbb D})$, the disc algebra.  It follows that
$X_{+} = X \cap A({\Bbb D})$, and so $X_{+} $ is a subalgebra of $X$.  Clearly
the maximal ideal space $M$ of $X_{+}$ contains clos ${\Bbb D}$, which is the maximal ideal space of $A({\Bbb D})$.  Since $X_{+}$ is generated as a Banach 
algebra by the single element $z$, $M$ is naturally identified with 
$\s_{X_{+}}(z)$, the spectrum of $z$ in $X_{+}$.  Since $X_{+}$ is a subalgebra
 of $X$ we have
$$
\partial \s_{X_{+}}(z) \subset \partial \s_X(z) = \partial \T = \T
$$
($\partial$ denotes boundary).  That is, $M$ contains clos ${\Bbb D}$ and
$\partial M \subset \T$.  Hence $M= {\rm clos}~ {\Bbb D}.\hfill \bl$

For a space $X$ of functions and a matrix-valued function $\Phi$ we write $\Phi \in X$
to mean that each entry of $\Phi $ belongs to $X$.  We denote by $X(M_{m,n})$ the space
of $m \times n$ matrix-valued functions whose entries belong to $X$, endowed with the
norm
$$
\| \Phi \|_X \df \sup \{ \|y^* \Phi x \|_X : \|x\|_{ \C^n} \le 1,~ \|y\|_{\C^m} \le 1 \}.
$$
$X(\C^n)$ is defined to be $X(M_{n,1})$.  For $\Phi \in L^\be(M_{m,n})$ we define the
Hankel operator $H_\Phi$ to be the operator from $H^2(\C^n)$ to $H^2_{-}(\C^m)$ given by
$$
H_\Phi x \df \pp_{-} (\Phi x).
$$
We recall that the space $QC$ of quasicontinuous functions is defined to be 
$(H^\be +C) \cap \overline{(H^\be +C)}$.

It transpires that the analysis of the continuity of $\A$ involves certain integer
indices associated with a matrix function.  These indices were introduced in [PY1], and
depend on the notion of a {\em thematic factorization}, which is a type of
diagonalization of a superoptimal error function $\Phi -\A \Phi$.  A {\em thematic function}
is a function $V \in  L^\be(M_{n,n})$ for some $n \in \N$ which is unitary-valued a.e.
on $\T$ and of the form
$$
V=\left(\begin{array}{cc}v &\bar{\a}\end{array}\right)
$$
where $v \in  H^\be(\C^n)$ is inner and co-outer 
and $\a \in  H^\be(M_{n,n-1})$ is co-outer.
Recall that an $H^\be$ matrix function $A$ is {\em inner} if $A(z)$ is an isometry for
almost all $z \in \T$ and is {\em co-outer} if $A^t H^2$ is dense in $H^2$.
Consider  $\Phi \in H^\infty+C$ of type $m \times n$.  We shall assume henceforth that
$m \le n$.  By [PY1, Theorem 2.1] the singular values $s_j(\Phi(z) - \A\Phi(z))$ are
constant a.e. on $\T$; their values $t_0 \ge t_1 \ge \dots \ge t_{m-1}$ are the {\em
superoptimal singular values} of $\Phi$.  Moreover, according to [PY1, Theorem 4.1],
$\Phi-\A\Phi$ admits a factorization of the form
\begin{equation}
\label{thfac}
\Phi-\A\Phi = W^*_0 W^*_1 \cdots W^*_{m-1} D V^*_{m-1} V^*_{m-2} \cdots V^*_0,
\end{equation}
where $D$ of type $m \times n$ is given by
$$
D \df \left(\begin{array}{ccccccc}t_0u_0 & 0 & \cdots & 0 & 0 & \cdots & 0 \\
0 & t_1u_1 & \cdots & 0 & 0 & \cdots & 0 \\
\cdot &\cdot &\cdots &\cdot &\cdot &\cdots & \cdot \\
0 & 0 &\cdots & t_{m-1}u_{m-1} & 0 & \cdots & 0
\end{array}\right)
$$
for some unimodular functions $u_0, \dots, u_{m-1} \in QC$,
$$
W_j =  \left(\begin{array}{cc} I_j & 0 \\ 0 & \tilde{W}_j \end{array}\right), \qquad
1 \le j \le m-1,
$$
$$
V_j =  \left(\begin{array}{cc} I_j & 0 \\ 0 & \tilde{V}_j \end{array}\right), \qquad
1 \le j \le m-1,
$$
and $W^t_0, ~ \tilde{W}_j^t,~V_0$ and $ \tilde{V}_j$ are thematic functions for 
$1 \le j \le m-1.$  We call (\ref{thfac}) a thematic factorization of $\Phi-\A\Phi$, and
we define the {\em index} of $t_j$ in this factorization to be the modulus of the
winding number of $u_j$ (or alternatively, as the Fredholm index of the Toeplitz
operator $T_{u_j}$).  Numerous properties of these indices were established in [PY3].
  In Section 1 we prove continuity of $\A$ with respect to a wide class of norms at
functions whose superoptimal singular values are nonzero and whose indices are all 1. 
For the Besov norm $B^1_1$ we obtain a continuity result even in the presence of zero
superoptimal singular values.  In Section 2 we consider the converse problem, and derive
some necessary conditions for continuity  points of $\A$
in the case of square matrix functions.  In
Section 3 we present sufficient conditions for the continuity of the superoptimal
singular values themselves.

\section{Sufficient conditions for continuity}
 
Let $X$ be a space of functions on $ \Bbb T$ invariant under $\A$ (e.g. one satisfying
the above axioms).
As we noted above, even in the scalar case $\A$
is discontinuous with respect to virtually any norm at any $\Phi$ such that $\| H_\Phi \|$
is a multiple singular value of $ H_\Phi $ [P1].  In the scalar case, for many spaces $X$
the converse also holds.  That is, if  $\| H_\Phi \|$ is a simple singular value then 
 $\Phi$ is a continuity point of $\A$ with respect to the norm of $X$.
 For matrix functions the situation is more complicated,
but we do have the following sufficient condition. 

\begin{thm}\thd
\label{t1.1}
Let $X$ be a space of functions on $\Bbb T$ satisfying Axioms {\rm (A1)} to {\rm (A5)}, 
let  $\Phi\in X(M_{m,n})$, $m\leq n$, and let $t_0,t_1,\cdots,t_{m-1}$ 
be the superoptimal singular values of $\Phi$. Suppose that $t_{m-1}\neq0$.
If $\Phi-{\cal A}\Phi$ has a thematic factorization with indices
\begin{equation}
\label{1.0}
k_0=k_1=\cdots=k_{m-1}=1,
\end{equation}
then $\Phi$ is a continuity point of the operator ${\cal A}$
of superoptimal approximation in $X(M_{m,n})$.
\end{thm}

As we have observed in [PY3], (\ref{1.0}) implies that all thematic factorizations of
$\Phi-\A \Phi$ have indices equal to 1.

The proof of the theorem will be based on the recursive construction of $\A
\Phi$ given in [PY2], which in turn was based on the proof in [PY1] that  $\A
\Phi$ is well defined.  Let us briefly recall the construction of  $\A \Phi$.
The first step is to find a Schmidt pair $\{v,~w\}$ of the compact Hankel
operator $H_\Phi$ corresponding to the singular value $\| H_\Phi \|$.   Then
find $Q \in H^\be (M_{m,n})$ such that
\begin{equation}
\label{qdef}
Qv = T_\Phi v, \qquad Q^t \bar{z} \bar{w} = T_{\Phi^t} \bar{z} \bar{w}
\end{equation}
(these equations always have a solution; in fact $Q= \A\Phi$ satisfies them, and in the
proof of the theorem below we shall even give an explicit rank two solution for $Q$).
Next let $v_{(i)}, ~ w_{(i)}$ be the inner factors of $v, ~  \bar{z}
\bar{w}$ and let
\begin{equation}
\label{VW}
V = \left(v_{(i)} \quad \bar{\alpha} \right), \qquad
 W^t = \left( w_{(i)} \quad \bar{\beta} \right)
\end{equation}
be thematic completions of  $v_{(i)}, ~ w_{(i)}$ respectively.  Then
$$
\A\Phi = Q + \beta \A \Psi \alpha^t
$$
where
$$
\Psi = \beta^* (\Phi - Q) \bar{\alpha}.
$$
Note that $\Psi$ is of type $(m-1) \times (n-1)$.  The strategy of the
proof is simply to show that $\a, \b$ and $Q$ can be chosen to depend
continuously on $\Phi$ and then to use induction on $m$.  In order to
do this we have to study some properties of maximizing
vectors for $H_\Phi$. 

It is easy to see from the axioms (A1)--(A5) that $H^*_\Phi H_\Phi$ 
is also a compact operator on $X_+(\C^n)$. 
Denote this operator on $X_+(\C^n)$ by $R$. 
We can identify the dual space $X_+^*(\C^n)$ with the space of 
analytic $\C^n$-valued functions $g$ in $\dd$ such that the Hermitian form
$$
(f,g)=\sum_{d\geq0}(\hat{f}(d),\hat{g}(d))_{\C^n}
$$
defined for polynomials $f$ in $X_+(\C^n)$, is continuous on $X_+(\C^n)$.
Obviously,
$$
X_+(\C^n)\subset H^2(\C^n)\subset X^*_+(\C^n).
$$
Since $R$ is a compact operator on $X_+(\C^n)$, it follows from the 
Riesz--Schauder theorem that $R^*$ is compact on $X^*_+(\C^n)$ and if $\l>0$, 
then $\l$ is an eigenvalue of $R$ if and only if $\l$ is an eigenvalue of $R^*$
of the same multiplicity (see [Yo], Ch. X, \S5). Since $H^*_\Phi H_\Phi$ is
self-adjoint, we have $R^*|H^2(\C^n)=H^*_\Phi H_\Phi$. Clearly, every
eigenvector of $R$ is an eigenvector of $H_\Phi^*H_\Phi$ and every eigenvector of
$H^*_\Phi H_\Phi$ is an eigenvector of $R^*$. It follows
 from the Riesz--Schauder theorem that $R$, $H^*_\Phi H_\Phi$, and $R^*$ have the same 
eigenvectors that correspond to positive eigenvalues.

\begin{thm}\thd
\label{t1.2}
Let $\Phi$ be a function in $X(M_{m,n})$, $m\leq n$, 
with superoptimal singular values
$t_0,\cdots,t_{m-1}$, $t_0\neq0$. Suppose that $\Phi-\A \Phi$ has
a thematic factorization whose indices $k_j$ are equal to 1 whenever 
$t_j=t_0$. Let $ \{v,~w \}$ be a Schmidt pair of $H_\Phi$ corresponding 
to $\|H_\Phi\|$. 
Then $v$ and $\bar{z} \bar{w}$ are co-outer and $v(\z)\neq0$ for any $\z\in\T$.
\end{thm}

Clearly it is sufficient to prove that $v(1)\neq\0$. 

We shall deduce Theorem \ref{t1.2} from the following lemma whose proof is similar to that of Lemma 3.2 of [PK].

\begin{lem}\thd
\label{t1.3}
Let $v$ be a maximizing vector for $H_\Phi$ such that $v(1)=0$.
Then $(1-z)^{-1}v\in X_+^*(\C^n)$ and $(1-z)^{-1}v$ is an eigenvector of $R^*$ with eigenvalue $t_0^2$.
\end{lem}

{\bf Proof.} Let us show that
\begin{equation}
\label{1.1}
(f,(1-z)^{-1}v)=(\pp_+v^*f)(1)
\end{equation}
for any polynomial $f$ in $X_+(\C^n)$, where
$v^*f(\z)\df(f(\z),v(\z))_{\C^n}$. Since the right-hand side of
(\ref{1.1}) clearly determines a continuous linear functional 
on $X_+(\C^n)$, it would follow that $(1-z)^{-1}v\in X_+^*(\C^n)$.

It is sufficient to establish (\ref{1.1}) for $f=z^jx$, $x\in\C^n$.
Obviously, $(1-z)^{-1}v=\sum_{j\geq0}z^j(\sum_{d=0}^j\hat{v}(d))$ and so
$$
(z^jx,(1-z)^{-1}v)=(x,\sum_{d=0}^j\hat{v}(d))_{\C^n}=\sum^j_{d=0}(x,\hat{v}(d))_{\C^n}.
$$
On the other hand it is easy to see that
$$
(\pp_+v^*f)(1)=\sum_{d=0}^j(x,\hat{v}(d))_{\C^n},
$$
which proves (\ref{1.1}).

To complete the proof of the lemma, 
we have to prove that $R^*(1-z)^{-1}v=t^2_0(1-z)^{-1}v$, which means that
\begin{equation}
\label{1.2}
(H^*_\Phi H_\Phi f,(1-z)^{-1}v)=t_0^2(f,(1-z)^{-1}v)
\end{equation}
for any $f\in X_+(\C^n)$. We may assume for convenience that $t_0=1$.

To establish (\ref{1.2}), we expand $(1-z)^{-1}v$ in the series
$\sum_{d=0}^\be z^dv$ and apply Ces\`{a}ro's summation method.

Let $Q\in H^\be$ be a best approximation to $\Phi$, i.e.
 $\|\Phi-Q\|_{L^\be}=\|H_\Phi\|=1$. Put $\Psi=\Phi-Q$.
 It is well known (see [AAK], [PY1], Th. 0.2) that
$H_\Phi v=\Psi v\in H^2_-(\C^m)$ and
$\|\Psi(\z)v(\z)\|_{\C^n}=\|v(\z)\|_{\C^n}$. Clearly, the last equality implies that $\Psi^*\Psi v=v$. 

We have
\bay
(H^*_\Phi H_\Phi f,z^dv)&=&(H^*_\Psi H_\Psi f,z^dv)=(H_\Psi f,\Psi z^dv)\nonumber\\
&=&(\pp_-\Psi f,z^d\Psi v)=(\Psi f,\pp_-z^d\Psi v)\nonumber\\
&=&(\Psi f,z^d\Psi v)-(\Psi f,\pp_+z^d\Psi v)\nonumber\\
&=&(f,z^d\Psi^*\Psi v)-(\Psi f,\pp_+z^d\Psi v)\nonumber\\
&=&(f,z^dv)-(\Psi f,\pp_+z^d\Psi v).\nonumber
\ey
Let $K_N(\z)\df\sum^N_{d=-N}(1-\frac{|d|}{N})\z^d$ be the Fej\'{e}r kernel and $K^+_N\df\pp_+K_N$. Then
$$
(H^*_\Psi H_\Psi f,K^+_Nv)=(f,K_N^+v)-(\Psi f,\pp_+K_N\Psi v),
$$
since $\Psi v\in H^2_-(\C^m)$. Let us prove that $\lim_{N\to\be}(\Psi f,\pp_+K_N\Psi v)=0$. Indeed
\bay
(\Psi f,\pp_+K_N\Psi v)&=&
(\pp_+\Psi f,K_N\Psi v)=(\Psi f,K_N\Psi v)-(\pp_-\Psi f,\Psi K_Nv)\nonumber\\
&=&(f,K_N\Psi^*\Psi v)-(H^*_\Psi H_\Psi f,K_Nv)\nonumber\\
&=&(f,K_Nv)-(H^*_\Psi H_\Psi f,K_Nv).\nonumber
\ey

Clearly,
$$
(f,K_Nv)\to(f(1),v(1))_{\C^n}=0;
$$
$$
(H^*_\Psi H_\Psi f,K_Nv)\to((H^*_\Psi H_\Psi f)(1),v(1))=0.
$$

It remains to prove that $\lim_{N\to\be}K_N^+v=(1-z)^{-1}v$ in the weak topology
$$\s(X_+^*(\C^n),X_+(\C^n)).$$ Let $g=z^jx$, $x\in\C^n$. Then
$$
(g,K_N^+v)=(v^*g,K_N^+)=
(\pp_+v^*g,K_N)\to(\pp_+v^*g)(1)
$$
as $N\to\be$.
The result follows now from (\ref{1.1}). \hfill $\bl$

\begin{cor}\thd
\label{t1.4}
Let $v$ be a maximizing vector for $H_\Phi$ such that $v(1)=\0$. Then
$(1-z)^{-1}v\in X_+(\C^n)$ and $(1-z)^{-1}v$ is also a maximizing vector for
$H_\Phi$.
\end{cor}

{\bf Proof of Theorem \ref{t1.2}.} Suppose that $v(1)=\0$. By Corollary 1.4, 
$v=(1-z)q$, where $q\in H^2$. Let
$$
w=\frac{1}{t_0}H_\Phi v.
$$
Then as we have already mentioned in the proof of Lemma \ref{t1.3},
$\|v(\z)\|_{\C^n}=\|w(\z)\|_{\C^n}$, $\z\in\T$. Let $h$ be a scalar outer
function such that $|h(\z)|=\|v(\z)\|$, $\z\in\T$, and let $h_1=(1-z)^{-1}h$.
Clearly, $h_1$ is also a scalar outer function and $|h_1(\z)|=\|q(\z)\|_{\C^n}$,
$\z\in\T$. Now there exist scalar inner functions $\vt_1$, $\vt_2$ such that
$v$ and $\bar{z}\bar{w}$ admit factorizations $v=\vt_1hv^{(i)}$, 
$\bar{z}\bar{w}=\vt_2hw^{(i)}$, where $v^{(i)}$ and $w^{(i)}$ are inner 
and co-outer in $H^2(\C^n)$. 
Then $hv^{(i)}$ is also a maximizing vector for $H_\Phi$ and 
$$
\frac{1}{t_0}H_\Phi hv^{(i)}=\bar{\vt}_1w
$$
(see the proof of Theorem 4.1 of [PY1]).

Let $V=\left(\begin{array}{cc}v^{(i)}&\bar{\a}\end{array}\right)$, 
$W^T=\left(\begin{array}{cc}w^{(i)}&\bar{\b}\end{array}\right)$ be thematic matrices. 
It follows from Lemma 2.3 of [PY1] that
$$
W(\Phi-\A\Phi)V=\left(\begin{array}{cc}t_0u_0&\0\\\0&F\end{array}\right)
$$
and $\Phi-\A\Phi$ has a thematic factorization with index equal to
$\dim\Ker T_{u_0}$, where 
\begin{equation}
\label{defu0}
u_0=\bar{z}\bar{\vt}_1\bar{\vt}_2\bar{h}/h.
\end{equation}
Since $h=(1-z)q$, we have
$$
u_0=\bar{z}\bar{\vt}_1\bar{\vt}_2\frac{1-\bar{z}}{1-z}\frac{\bar{q}}{q}
=-\bar{z}^2\bar{\vt}_1\bar{\vt}_2\frac{\bar{q}}{q},
$$
and so $k_0=\dim\Ker T_{u_0}\geq2$, since obviously $q$ and $zq$ belong to
$\Ker T_{u_0}$. This contradicts the hypotheses of Theorem 1.2, and so
$v(1) \ne 0$.  In similar fashion, the relation (\ref{defu0}) shows that
 Ker $T_{u_0}$ contains $h, ~\bar{\vt_1} h$ and $\bar{\vt_2} h$.  Thus,
 if $v^{(i)}$ or $w^{(i)}$ is not co-outer, we have again contradicted $\dim~ {\rm Ker}
~T_{u_0}=1$.  Hence $v,~w$ are co-outer. $\hfill \bl$

\begin{lem}\thd
Let $n>1$ and let $\f$ be an inner function in $X_+(\C^n)$.
Then $0$ is an isolated spectral point of the operators
$T^X_{\bar\f}T^X_{\f^t}$ on $X_+(\C^n)$ and $T_{\bar\f}T_{\f^t}$ on $H^2(\C^n)$.
\end{lem}

\pf Let us prove the lemma for the operator $T^X_{\bar\f}T^X_{\f^t}$.
The proof for $T_{\bar\f}T_{\f^t}$ is exactly the same. 

Let us observe that we may assume that $\f$ is co-outer. Indeed if
$\f=\vt\t$, where $\vt$ is a scalar inner function and $\t$ is an inner co-outer
function, then it follows from the axiom (A5) that $\t\in X_+(\C^n)$ and
clearly $T^X_{\bar\f}T^X_{\f^t}=T^X_{\bar\t}T^X_{\t^t}$.

Consider the operator
$T^X_{\f^t}T^X_{\bar\f}$ on $X_+$. It is well known that a nonzero point $\l\in\C$
belongs to the spectrum of $T^X_{\bar\f}T^X_{\f^t}$ if and only if it belongs
to the spectrum of $T^X_{\f^t}T^X_{\bar\f}$. Therefore to prove the lemma it is 
sufficient to show that $T_{\f^t}^XT_{\bar\f}^X$ is invertible.

We have 
$$
T_{\f^t}^XT_{\bar\f}^X=I-H^{*X}_{\bar\f}H^X_{\bar\f}.
$$
It follows easily from the axioms (A1)--(A5) that the operator 
$H^{*X}_{\bar\f}H^X_{\bar\f}$ is compact. Hence it is sufficient to show that
$\Ker T_{\f^t}^XT_{\bar\f}^X=\bO$. Let $f\in\Ker T_{\f^t}^XT_{\bar\f}^X$.
Then $H^{*X}_{\bar\f}H^X_{\bar\f}f=f$, which clearly means that
$H^*_{\bar\f}H_{\bar\f}f=f$. Since 
$\|H_{\bar\f}\|=1$ and $\|\bar\f\|_{L^\be(\C^n)}=1$, it follows
that $\bar\f f\in H^2_-(\C^n)$.
Thus $\bar{f} \f^t H^2(\C^n) \subset zH^1$, and since $\f^tH^2(\C^n)$
is dense in $H^2$, it follows that $\bar f H^2 \subset zH^1$, and hence that
 $\bar f \in zH^2$.   Thus $f=\0$. $\hfill \bl$

For an inner function $\f\in H^\be(\C^n)$ we denote by $L_\f$ the
kernel of $T_{\f^t}$ and by $P_\f$ the orthogonal projection from
$H^2(\C^n)$ onto $L_\f$. 
Similarly, we denote by $L_\f^X$ the kernel
of $T^X_{\f^t}$. Clearly, $L_\f=\Ker T_{\bar\f}T_{\f^t}$ and
$L_\f^X=\Ker T^X_{\bar\f}T^X_{\f^t}$.

Consider a simple closed  positively oriented Jordan curve $\O$ which lies in
the resolvent sets of $T_{\bar\f}T_{\f^t}$ and $T^X_{\bar\f}T^X_{\f^t}$,
encircles zero but does not wind round any other point of the spectra of
 $T_{\bar\f}T_{\f^t}$ and $T^X_{\bar\f}T^X_{\f^t}$. Clearly
$$
P_\f=\frac{1}{2\pi  \rm i}\oint_\O(\z I - T_{\bar\f}T_{\f^t})^{-1}d\z.
$$
Consider the projection $P^X_\f$ from $X_+(\C^n)$ onto $L_\f^X$ defined by
\begin{equation}
\label{1.4}
P^X_\f=\frac{1}{2\pi \rm i}\oint_\O(\z I - T^X_{\bar\f}T^X_{\f^t})^{-1}d\z.
\end{equation}
Obviously, $P^X_\f f=P_\f f$ for $f\in X_+(\C^n)$.

Suppose now that $\{\f^{(k)}\}_{k\ge1}$ is a sequence of inner functions in
$X_+(\C^n)$, which converges to $\f$ in the norm.
Then  $T^X_{(\f^{(k)})^t}T^X_{\bar\f^{(k)}} \to  T^X_{\f^t}T^X_{\bar\f}$ in the norm
of ${\cal L}(X_{+}(\C^n)).$ As in the proof of Lemma 1.5, $ T^X_{\f^t}T^X_{\bar\f}$
is invertible, and hence there is a neighbourhood $U$ of zero which lies in the
resolvent set of $ T^X_{\f^t}T^X_{\bar\f}$ and of  
$T^X_{(\f^{(k)})^t}T^X_{\bar\f^{(k)}}$ for all sufficiently
large $k$.  Without loss of generality we may assume that this
holds for all values of $k$.  Choose a simple closed contour $\O$
lying in $U$ and winding round $0$.
Then $0$ is the only point inside or on $\O$ of the spectra of 
$T_{\bar\f^{(k)}}T_{(\f^{(k)})^t}$ and
$T^X_{\bar\f^{(k)}}T^X_{(\f^{(k)})^t}$.
We can therefore define projections $P_{\f}, ~P^X_\f,~P_{\f^{(k)}}, ~P^X_{\f^{(k)}}$ 
by integrals as above, all using the same contour $\O$.
It is then easy to see from (\ref{1.4}) that
$P^X_{\f^{(k)}}\to P^X_\f$ in the operator norm. 

\begin{lem}\thd
Let $V=\left(\begin{array}{cc}\f&\ov{\f_c}\end{array}\right)$ be
unitary-valued on $\T$, where
$\f_c$ is inner and co-outer. There exist inner co-outer functions  $\f_c^{(k)}$
such that
$V^{(k)} \df \left(\begin{array}{cc}\f^{(k)}&\ov{\f^{(k)}_c}\end{array}\right)$
is unitary-valued on $\T$ and $\|V-V^{(k)}\|_{X(M_{n,n})}\to0$.
\end{lem}

\pf It was shown in [PY1] (see the proof of Theorem 1.1) that, for a given
inner column $\f$, one can construct an inner co-outer $\a$ such that 
$\left(\begin{array}{cc}\f&\ov{\a}\end{array}\right)$ is unitary-valued on $\T$ and
the columns of $\a$ have the form
$P_\f C_1,P_\f C_2,\cdots,P_\f C_{n-1}$, where 
$C_1,C_2,\cdots,C_{n-1}$ are constant column functions.
By [PY1, Corollary 1.6], $\f_c = \a U$ for some constant unitary $U$.  Hence
the columns of $\f_c$ also have the form $P_\f C_j$ for some constants $C_j$.
Consider the subspace of $H^2(\C^n)$
$$
P_\varphi \C^n  \df \{P_\varphi C : C \in \C^n \},
$$
where we identify $C \in \C^n$ with a constant function in $H^2(\C^n)$.
This space has the remarkable property that the pointwise and $H^2$ inner
products coincide on it.  That is, if $f_j = P_\varphi C_j,~ j=1,2,$ where
$C_1, ~C_2 \in \C^n$, then
\begin{equation}
\label{ptwh2}
(f_1,f_2)_{H^2(\C^n)} = (f_1(z),f_2(z))_{\C^n}
\end{equation}
for almost all $z \in \T$.  To see this note that $L_\varphi$ is a closed
$z$-invariant subspace of $H^2(\C^n)$, and so is of the form $\Theta
H^2(\C^p)$ for some natural number $p$ and some $n \times p$ inner function
$\Theta$. Then for any $C \in  \C^n$,
$$
P_\varphi C = \Theta \pp_{+} \Theta^* C = \Theta \Theta(0)^* C,
$$
and so
\bay
(f_1,f_2)_{H^2(\C^n)}& =& (P_\varphi C_1, P_\varphi C_2)_{H^2}
 = ( \Theta \Theta(0)^* C_1,  \Theta \Theta(0)^* C_2)_{H^2}  \nonumber \\
 &=& (\Theta(0)^* C_1,\Theta(0)^* C_2)_{H^2} =
 (\Theta(0)^* C_1,\Theta(0)^* C_2)_{\C^p}  \nonumber \\
&=& ( \Theta(z) \Theta(0)^* C_1,  \Theta(z) \Theta(0)^* C_2)_{\C^n}
  \nonumber \\
&=&  (f_1(z),f_2(z))_{\C^n}            \nonumber
\ey
for almost all $z \in \T$.
It follows that any unit vector in $P_\varphi \C^n$ is an inner column
 function, and any orthonormal sequence (with
respect to the inner product of $H^2 (\C^n)$) of vectors in  $P_\varphi
\C^n$ constitutes the columns of an inner function.

Now let $P_\varphi C_j,~ 1 \le j \le n-1$, be the columns of $\f_c$ as above,
and consider the functions
$P_{\f^{(k)}}C_1,P_{\f^{(k)}}C_2,\cdots,P_{\f^{(k)}}C_{n-1}$. 
Clearly
$$
\|P_\f C_j-P_{\f^{(k)}}C_j\|_{X(\C^n)}=
\|P^X_\f C_j-P^X_{\f^{(k)}}C_j\|_{X(\C^n)}\to0\quad\mbox{as}\quad k\to\be.
$$
It follows that for large values of $k$ the inner products
$(P_{\f^{(k)}} C_{j_1},P_{\f^{(k)}} C_{j_2})_{H^2(\C^n)}$
 are small for $j_1\neq j_2$ and are
close to 1 if $j_1=j_2$.   We shall show that the desired $\f_c^{(k)}$ can
be obtained by orthonormalising the  $ P_{\f^{(k)}}C_j$.

Pick $M > 1$ such that $\| P_{\f^{(k)}}C_j\|_{X(\C^n)} \le M$ for all $k \in
\N$ and $1 \le j < n$.  By the equivalence of norms on finite-dimensional
spaces there exists $K > 0 $ such that, for any $(n-1)$-square matrix
$T = \left( t_{ij} \right)$,
\begin{equation}
\label{norms}
\max \mid  t_{ij} \mid ~ \le ~ \| T \| ~ \le  ~K \max \mid  t_{ij} \mid
\end{equation}
(here $\|.\|$ is the operator norm on ${\cal L} ( \C^{n-1})$).

Let $0 < \e < 1$.  Choose $k_0$ such that $k \ge k_0$ implies
\begin{equation}
\label{cols}
\| P_{\f^{(k)}}C_j - P_\varphi C_j \|_{X(\C^n)} < \frac{\e}{2}, \qquad
j=1,\dots,n-1,
\end{equation}
and
\begin{equation}
\label{entries}
\mid ( P_{\f^{(k)}}C_i,~ P_{\f^{(k)}}C_j ) - \d_{ij} \mid < \frac{\e}{2KnM},
\qquad i,j= 1,\dots,n-1.
\end{equation}
Fix $k \ge k_0$ and let $T: \C^{n-1} \rightarrow  P_{\f^{(k)}}\C^n$ be the
operator which maps the $j$th standard basis vector $e_j$ of $\C^{n-1}$ to
  $ P_{\f^{(k)}}C_j$.  The matrix of $T^*T \in {\cal L} ( \C^{n-1})$ is the
Gram matrix $ ( P_{\f^{(k)}}C_j,~ P_{\f^{(k)}}C_i )$, and so by
(\ref{norms}) and (\ref{entries}) we have
$$
\|T^*T - I \| < \frac{\e}{2nM} < \frac{1}{2}.
$$
By diagonalisation,
$$
\|(T^*T)^{-\frac{1}{2}} - I \| < \frac{\e}{2nM}.
$$
Let $(T^*T)^{-\frac{1}{2}} = \left( t_{ij} \right)$: then
$$
\mid  t_{ij} -  \d_{ij} \mid <  \frac{\e}{2nM}.
$$
Let the polar decomposition of $T$ be $T=U(T^*T)^{\frac{1}{2}}$, so that 
$U=T(T^*T)^{-\frac{1}{2}}.$   Then $U$ is unitary, so that
$Ue_1,\dots,Ue_{n-1}$ are orthonormal in $ P_{\f^{(k)}}\C^n$.  Let
 $\f_c^{(k)}$ be the $n \times (n-1)$ matrix with columns
$Ue_1,\dots,Ue_{n-1}.$  By the remark above,  $\f_c^{(k)}$ is inner.  By the
fact that
 $ P_{\f^{(k)}}\C^n \subset L_{\f^{(k)}}$, the columns of
$\bar{\f}_c^{(k)}$
are pointwise orthogonal to $\f^{(k)}$. Hence 
$$
V^{(k)} \df \left( \f^{(k)} ~\bar{\f}_c^{(k)} \right)
$$
is unitary-valued. 
Furthermore, the $j$th column $Ue_j$ of $\f_c^{(k)}$ satisfies
\bay
\| P_{\f^{(k)}} C_j - Ue_j \|_{X(\C^n)}& =& \|Te_j - T(T^*T)^{-\frac{1}{2}}e_j
\|_X \\ & = &
\|Te_j -\left(Te_1 \dots Te_{n-1} \right) \left( \begin{array}{c} t_{1j} \nonumber  \\
\vdots  \\ t_{n-1,j} \end{array} \right) \|_X  \nonumber \\
&\le& \mid t_{1j} \mid \|Te_1\|_X +\dots + \mid t_{jj} -1 \mid \|Te_j \|_X + \nonumber  \\
 & & \hspace{.3cm} \dots + \mid t_{n-1,j} \mid \|Te_{n-1} \|_X  \nonumber \\
& \le &  (n-1) \frac{\e}{2nM} M <  \frac{\e}{2}. \nonumber
\ey
On combining this inequality with (\ref{cols}) we obtain
\bay
\| P_\f C_j - Ue_j \|_X  & \le & \|  P_\f C_j -  P_{\f^{(k)}} C_j \|_X +
\|  P_{\f^{(k)}} C_j - Ue_j \|_X \nonumber \\
& \le & \frac{\e}{2} + \frac{\e}{2} = \e. \nonumber
\ey
That is, the $j$th column of  $\f_c^{(k)}$ tends to the $j$th column of
$\f_c$ with respect to the norm of $X(\C^n)$. Hence
$V^{(k)} \rightarrow V$ in $X(M_{n,n})$.
Finally, it follows from  [P3, Lemma 1.2] that  $\f_c^{(k)}$ is co-outer.  $ \hfill \bl$

\begin{cor}\thd
Suppose $\f$ is co-outer and $V^{(k)}$ is constructed as in Lemma 1.6.  For
sufficiently large $k$,  $\f^{(k)}$ is co-outer and so $V^{(k)}$ is thematic.
\end{cor}

\pf By [PY1, Theorem 1.2], $\det V$ is constant, hence has zero winding 
number about 0.  Since $\det V^{(k)} \to \det V$ uniformly on $\T$,
 $\det V^{(k)}$ also has zero winding number about 0 for sufficiently 
large $k$.  Again by [PY1, Theorem 1.2], $\f^{(k)}$ is co-outer. $ \hfill \bl$

\begin{lem}\thd
Let $E,~F$ be Banach spaces, let $T:E \ra F$ be a surjective continuous
 linear mapping
and let $x \in E,~y \in F$ be such that $Tx=y$.  Let $\e > 0$ and let $T'
\in {\cal L} (E,F)$.  There exists $\d > 0$ such that, whenever $\| T' -T \|
< \d$, the equation $T'x'=y$ has a solution $x'$ satisfying $\|x' - x\| < \e$.
\end{lem}

\pf  We can suppose $\|x\| = 1$.  By the Open Mapping Theorem there exists
$c > 0$ such that the ball of radius $c$ in $F$ is contained in the image
under $T$ of the unit ball in $E$.  Then
$$
\|T^*f \| \ge c \|f\| \qquad {\rm for ~ all ~} f \in F^*.
$$
Let $\d = \frac{c}{2} \min \{1,\e \}.$  Suppose  $\| T' -T \| < \d$.
  For any $f \in F^*$
\bay
\|T'{}^*f \| &=& \|T^*f + (T'-T)^*f \| \ge \|T^*f\| - \|T' - T\| \cdot \|f\|
\nonumber \\
&\ge& c\|f\| -  \frac{c}{2} \|f\| = \frac{c}{2} \|f\|.
\ey
Thus $T'$ maps the closed unit ball of $E$ to a superset of the closed ball
of radius $\frac{c}{2}$ in $F$.  Since
$$
\|(T -T')x \| \le \|T - T' \| < \d,
$$
it follows that there exists $\xi \in E$ such that
$$
\| \xi \| < \frac{2\d}{c} \le \e
$$
and $T' \xi = (T - T')x$.  Then $x' \df x + \xi $ has the stated properties:
$$
T'x' = T'x + (T - T')x = Tx = y,
$$
\hspace{5.5cm} $
\|x - x' \| = \| \xi \| < \e. \hfill  \bl$

\begin{lem}\thd
Let $f,~\f \in X_{+}(\C^n)$ be such that $\f^tf = 1$ and let $\e > 0$.  
There exists $\d > 0$ such that, for any $\psi \in  X_{+}(\C^n)$ satisfying
$\|\f - \psi \|_X < \d$, there is a $g \in  X_{+}(\C^n)$ such that
$\|f-g\|_X < \e$ and $\psi^t g = 1$.
\end{lem}

\pf  Let $T = T^X_{\f^t} : X_{+}(\C^n) \ra X_{+}$, so that $Tx = \f^t x$ for
$x \in X_{+}(\C^n)$.  Then $T$ is a surjective continuous linear mapping and
$Tf = 1$.  By Lemma 1.7 there exists $\d$ such that $\|T' - T\| < \d$
implies that the equation $T'g = 1$ has a solution $g \in X_{+}(\C^n)$
satisfying $\|f - g\|_X < \e$.   If $\psi \in X_{+}(\C^n)$ is such that
$\|\f - \psi\|_X < \d$ then $\|T^X_{\psi^t} - T^X_{\f^t} \| < \d$, and so
the lemma applies to $T' = T^X_{\psi^t}$; that is, there exists $g \in 
 X_{+}(\C^n)$ such that $\psi^t g = 1$ and $\|f - g \|_X < \e. \hfill \bl$

{\bf Proof of Theorem \ref{t1.1}.} We proceed by induction on $m$. 
% Suppose first that $m>1$.

Let $\{\Phi^{(k)}\}_{k\geq1}$ be a sequence of functions in $X$ such that 
$\|\Phi-\Phi^{(k)}\|_{X(M_{m,n})}\rightarrow0$.  We shall show that some
subsequence of $\A \Phi^{(k)}$ converges to $\A \Phi$ in the norm of $X$:
this will suffice to establish the continuity of $\A$ at $\Phi$.
Let $v^{(k)}$ be a co-outer maximizing vector for the operator 
$H_{\Phi^{(k)}}$ on $H^2(\C^n)$.
We can take it that the norm of  $v^{(k)}$ in $X_+(\C^n)$  is equal to 1:
$$
\|v^{(k)}\|_{X(\C^n)}=1,~~~~\|H_{\Phi^{(k)}}v^{(k)}\|_{H_-^2(\C^m)}=
\|H_{\Phi^{(k)}}\|\cdot\|v^{(k)}\|_{H^2(\C^n)}.
$$

Let $\O$ be a  positively oriented Jordan contour which
winds once round the largest eigenvalue $t^2_0$ of $H_\Phi^*H_\Phi$, 
contains no eigenvalues and encircles no other
eigenvalues. 

It is easy to see from the axioms $(A1)$--$(A5)$ that the operators
$H_{\Phi^{(k)}}^{*X}H^X_{\Phi^{(k)}}$ converge to 
$H_\Phi^{*X}H^X_\Phi$  in the operator norm of $X_+(\C^n)$. 
It follows that for large values of $k$ there are no points of the spectrum
of $H_{\Phi^{(k)}}^{*X}H^X_{\Phi^{(k)}}$ on $\O$. Let
$$
\p=\frac{1}{2\pi \rm i}\oint_\O(\z I-H_\Phi^{*X}H^X_\Phi)^{-1}d\z
$$
and
$$
\p^{(k)}=
\frac{1}{2\pi \rm i}\oint_\O(\z I - H_{\Phi^{(k)}}^{*X}H^X_{\Phi^{(k)}})^{-1}d\z.
$$

Clearly, $\p v^{(k)}$ is a maximizing vector of $H_\Phi^{*X}H^X_\Phi$ and
$$
\|v^{(k)}-\p v^{(k)}\|_{X(\C^n)}=\|\p^{(k)}v^{(k)}-\p v^{(k)}\|_{X(\C^n)}\to0,
\quad k\to\be.
$$

The vectors $\p v^{(k)}$ belong to the finite-dimensional subspace of
maximizing vectors of $H_\Phi^*H_\Phi$. Therefore there exists a convergent
subsequence of the sequence $\{\p v^{(k)}\}_{k\ge0}$. Without loss of generality
we may assume that the sequence $\{\p v^{(k)}\}_{k\ge0}$ converges in
$X(\C^n)$ to a vector $v$, which is a maximizing vector of $H_\Phi^*H_\Phi$.
Obviously, $\|v^{(k)}-v\|_{X(\C^n)}\to0$ as $k\to\be$.

We also need the other Schmidt vectors corresponding to $v$ and $v^{(k)}$.
We may assume that $ \| H_{\Phi^{(k)}} \| \ne 0$ for all $k$.  Let
$$
w = t_0^{-1} H_\Phi v, \qquad w^{(k)} =  H_{\Phi^{(k)}} v^{(k)} /
 \| H_{\Phi^{(k)}} \|
$$
in $X_{-}(\C^m).$
Since $ H_{\Phi^{(k)}} \to  H_\Phi$ in the norm of $ {\cal L}(X_{+}(\C^n),
X_{-}(\C^m))$ and $v^{(k)} \to v$ in $X_{+}$ it follows that $w^{(k)} \to w$
in $X_{-}$.  The $v^{(k)}$ are co-outer by choice; the same is true of $w^{(k)}$
for sufficiently large $k$ by Corollary 1.7.

Now let us show that Theorem \ref{t1.1} holds when $m=1$.  In this case
$w$ and $w^{(k)}$ are scalar functions in $X$.  By [AAK],
 $ |w(z)| = \| v(z) \|$ a.e. on $\T$.  By continuity, equality holds at
all points of $\T$.  By Theorem \ref{t1.2}, $v$ (and hence also $w$) is
non-zero at every point of the maximal ideal space $\T$ of $X$.  Thus $1/w
\in X$.  By virtue of the continuity of inversion in Banach algebras we
deduce that $1/w^{(k)} \in X$ for sufficiently large $k$, and  $1/w^{(k)}
\to 1/w$ in $X$.  Again by [AAK],
$$
w^*(\Phi - \A \Phi) = \| H_\Phi\| v^*
$$
and hence
$$
\Phi - \A \Phi =  \| H_\Phi\| \frac{v^*}{w^*}, \qquad
\Phi^{(k)} - \A \Phi^{(k)} =  \| H_{\Phi^{(k)}}\| \frac{v^{(k)*}}{w^{(k)*}},
$$
the latter for large $k$.  From these equations it is clear that 
$ \A \Phi^{(k)} \to  \A \Phi$ in $X$.  Thus the case $m=1$ is established.

Now consider $m>1$ and suppose the theorem true for $m-1$.  We prove the
induction step by block-diagonalisation of $\Phi - \A \Phi$. 
Let $v,~w$ be as above and let $h$ be the outer factor of $v$.  Once again by [AAK],
$h$ is also the outer factor of $\bar{z} \bar{w}$. 
 It is given explicitly by the formula [H]
$$
h = e^{u + {\rm i} \tilde{u}}
$$
where
$$
u= \log \|v(\cdot) \|
$$
and $\tilde{u}$ is the harmonic conjugate of $u$,
$$
\tilde{u} = - {\rm i} (2\pp_{+} - I)u.
$$
Since $v \in X_{+}(\C^n)$ it is clear from axioms
 A1 and A2 that $\|v(\cdot)\|^2 \in X$.  By
Theorem 1.2,  $\|v(\cdot)\|^2$ does not vanish on $\T$, and so its spectrum in the Banach
algebra $X$ is a compact interval of the positive real numbers.  By the analytic
functional calculus, $u = \frac{1}{2} \log  \|v(\cdot)\|^2 \in X$.  By A1 we have also
$\tilde{u} \in X$.
Thus $h = e^{u + {\rm i} \tilde{u}} \in X$.
The above construction also makes it clear that if $v^{(k)},~ h^{(k)}$ are the
corresponding entities for $\Phi^{(k)}$, so that $v^{(k)} \rightarrow v$ in $X$, then
$ h^{(k)} \ra h$ in $X$.  Indeed, since $\pp_{+}$ maps $X$ into itself, it follows from
the Closed Graph Theorem that  $\pp_{+}$ is continuous on $X$, and hence the Hilbert 
transform $u \ra \tilde{u}$ is continuous on $X$.

Note also that since $ |h| = \|v(\cdot) \|$ is bounded away from zero, $h$ is invertible in
$X$ and $1/h^{(k)} \ra 1/h$ in $X$.  Let 
$v_{(i)},~ w_{(i)}, ~v_{(i)}^{(k)},~ w_{(i)}^{(k)}$ be the
inner factors of $v,~ \bar{z} \bar{w},~ v^{(k)}, ~ \bar{z} \bar{ w}^{(k)}$
respectively, so that
$$
v_{(i)} = v/h, \qquad  v_{(i)}^{(k)} =  v^{(k)}/ h^{(k)}
$$
etc. Then
$  v_{(i)}^{(k)} \ra v_{(i)}$ and $ w_{(i)}^{(k)} \ra w_{(i)}$ in $X_{+}$
as $k \ra \be$. By Theorem \ref{t1.2}, $v_{(i)}$ and $w_{(i)}$ are co-outer.

By Lemma 1.6 we can find thematic functions
$$
V = \left(v_{(i)} \quad \bar{\a} \right), \qquad W^t = \left(  w_{(i)} \quad
\bar{\b} \right),
$$
$$
V^{(k)} = \left( v_{(i)}^{(k)} \quad \bar{\a}^{(k)} \right), \qquad
W^{(k)t} =  \left( w_{(i)}^{(k)} \quad \bar{\b}^{(k)} \right)
$$
such that $V^{(k)} \to V$ and $W^{(k)} \to W$ in $X$. {\em  A fortiori},
\begin{equation}
\label{albet}
\a^{(k)} \to \a, \qquad \b^{(k)} \to \b
\end{equation}
in $X(M_{n,n-1}), ~ X(M_{m,m-1})$ respectively.

Now we construct $Q,~ Q^{(k)} \in  X_{+}(M_{m,n})$ such that (cf
(\ref{qdef}))
\begin{equation}
\label{qdef2}
Qv = T_\Phi v, \qquad Q^t\bar{z}\bar{w} = T_{\Phi^t} (\bar{z}\bar{w}),
\end{equation}
\begin{equation}
\label{qdef3}
Q^{(k)}v^{(k)} = T_{\Phi^{(k)}} v^{(k)}, \qquad
 Q^{(k)t}\bar{z}\bar{w}^{(k)} = T_{\Phi^{(k)t}} (\bar{z}\bar{w}^{(k)})
\end{equation}
and $Q^{(k)} \to Q$ in $X$.   We can do this using a formula for $Q$ which
we gave in [PY2, Sec. 2, Remark 3].  Let
$$
y_1 = T_\Phi v/ h, \qquad y_2 =  T_{\Phi^t} (\bar{z} \bar{w}) /h.
$$
Then $y_1,~y_2 \in X$ and from the fact that the equations (\ref{qdef2}) are consistent
 (they hold with $Q = {\cal A}\Phi$) we have
$y_2^t v_{(i)} = w_{(i)}^t y_1 ~(=  w_{(i)}^t Q  v_{(i)})$.
The components of $v_{(i)}$ are elements of the Banach algebra $X_{+}$.  By
Theorem 1.2 they do not vanish simultaneously at any point of $\T$, nor (since
$v_{(i)}$ is co-outer) do they at any point of $\Bbb D$.  Hence they 
do not all 
belong to any maximal ideal of $X_{+}$ (see Lemma 0.1), and so the ideal 
they  generate in $X_{+}$ is the whole algebra.  Thus there exists
$f_1 \in X_{+}(\C^n)$ such that $f_1^t v_{(i)} =1$.  Likewise there exists
$f_2 \in X_{+}(\C^m)$ such that $f_2^t w_{(i)} =1.$
%$$
%f_1^t  v_{(i)} = 1 = f_2^t  w_{(i)}.
%$$
It is simple to verify that a solution of (\ref{qdef2}) is
$$
Q = y_1 f_1^T + f_2 y_2^t - f_2 y_2^t v_{(i)} f_1^t.
$$
Now perform a similar construction to obtain $Q^{(k)}$.
Let
$$
y_1^{(k)} = T_{\Phi^{(k)}} v^{(k)}/ h^{(k)}, \qquad 
y_2^{(k)} =  T_{\Phi^{(k)t}} (\bar{z} \bar{w}^{(k)}) /h^{(k)}.
$$
Then $y_1^{(k)} \to y_1$ and $y_2^{(k)} \to y_2$ in $X$.

Apply Lemma 1.8 to $f=f_1, ~ \f = v_{(i)}$.  For any $N \in \N$ there exists
$\d_N > 0$ such that $\| v_{(i)} - \psi \|_X < \d_N$ implies that there
exists $g \in X_{+}(\C^n)$ with $g^t\psi = 1$ and $\|f_1 - g\|_X <
\frac{1}{N}$.  Define a sequence of integers $(k_N)$ and $f_1^{(k_N)} \in
X_{+}(\C^n)$ inductively as follows.   Let $k_1 =1, \quad f_1^{(1)} =0$.  Choose
$k_N > k_{N-1}$ so that $\|v_{(i)} - v_{(i)}^{(k_N)} \| < \d_N$.  Then there
exists $f_1^{(k_N)} \in X_{+}(\C^n)$ such that $f_1^{(k_N) t}  v_{(i)}^{(k_N)} =1$ and
$\|f_1^{(k_N)} - f_1\| < \frac{1}{N}$.   Passing to the subsequence
$(\Phi^{(k_N)})$ of $(\Phi^{(k)})$, we may assume that $f_1^{(k)t}
v_{(i)}^{(k)} = 1$ and $f_1^{(k)} \to f_1$ in $X$.  In a similar way we
construct $f_2^{(k)} \in X_{+}(\C^m)$ such that  $f_2^{(k)t}
w_{(i)}^{(k)} = 1$ and $f_2^{(k)} \to f_2$ in $X$.   Now let
$$
Q^{(k)} = y_1^{(k)} f_1^{(k)t} + f_2^{(k)} y_2^{(k)t}
 - f_2^{(k)} y_2^{(k)t} v_{(i)}^{(k)} f_1^{(k)t}.
$$
Then $Q^{(k)}$ satisfies (\ref{qdef3}) and $Q^{(k)} \to Q$ in $X$.  Let
$$
\Psi \df \b^*(\Phi - Q) \bar{\a}, \qquad
\Psi^{(k)} \df \b^{(k)*}(\Phi^{(k)} - Q^{(k)}) \bar{\a}^{(k)}.
$$
Then $\Psi^{(k)} \to \Psi$ in $X(M_{m-1,n-1})$.  It is shown in [PY1,PY2]
that
\begin{equation}
\label{recur1}
\Phi - \A \Phi = W^* \left( \begin{array}{cc} t_0 u_0 & 0 \\ 0 & 
\Psi - \A \Psi \end{array} \right) V^*
\end{equation}
where $u_0$ is a badly approximable unimodular function.  It follows
that the superoptimal singular values of $\Psi$ are $t_1, \dots, t_{m-1}$
and are non-zero.  Furthermore, every thematic factorization of $\Psi - \A
\Psi $ gives rise to one of $\Phi - \A \Phi$, and hence the indices in any
 thematic factorization of $\Psi - \A
\Psi $ are all equal to 1.   By the inductive hypothesis $\A$ is
continuous at $\Psi$, and hence $\A \Psi^{(k)} \to \A\Psi$ in $X(M_{m-1,n-1})$.
By [PY2],
\begin{equation}
\label{recur2}
\A \Phi = Q + \b \A\Psi \a^t, \qquad
\A \Phi^{(k)} = Q^{(k)} + \b^{(k)} \A\Psi^{(k)} \a^{(k)t}
\end{equation}
and hence $\A \Phi^{(k)} \to \A \Phi$ in $X(M_{m,n})$ as $k \to \be$. Thus
$\A$ is continuous at $\Phi. \hfill \bl$

What if one of the superoptimal singular values $t_j$ of $\Phi$ is 0?  One can see by
considering diagonal examples such as diag$\{\bar{z},0 \}$ that it is important whether
$\A$ is continuous at 0 in (scalar) $X$, or equivalently whether $\A$ is bounded.
  This is not always so for spaces
satisfying A1 to A5 (see [P2]), and so the conclusion of Theorem 1.1 does not follow 
if the condition $t_{m-1} \neq 0$ is relaxed. There is one case when it does.

\begin{thm}\thd
\label{t1.9}
Let $X$ be the Besov space $B_1^1$ and let
 $\Phi\in X(M_{m,n})$.
If \mbox{ $\Phi-\A \Phi$} has a thematic factorization in which the indices corresponding
to non-zero superoptimal singular values are all equal to 1
then $\Phi$ is a continuity point of the operator ${\cal A}$
of superoptimal approximation in $X(M_{m,n})$.
\end{thm}

\pf
 The fact that this statement is true in the case $\Phi = \0$ is Theorem 5.6 of [PY1].
 Note that $X$ satisfies axioms A1 to A5.  Let  $\Phi$ have
superoptimal singular values $t_0,\dots,t_{m-1}$.
Let $r$ be the number of nonzero superoptimal singular values of $\Phi$:  
$r = \inf \{j:t_j = 0\}.$  We prove the result by induction on $r$.
As in the proof of Theorem \ref{t1.1},
let $\{\Phi^{(k)}\}_{k\geq1}$ be a sequence of functions in $X$ such that 
$\|\Phi-\Phi^{(k)}\|_{X(M_{m,n})}\rightarrow0$.

  If $r=0$ then $\Phi =\A \Phi \in H^\infty$.  Since $\Phi^{(k)} \rightarrow \Phi$ in
$X$, by the cited theorem,  $\A (\Phi^{(k)} - \Phi) \rightarrow  \0$, and since $\Phi \in
 H^\infty$,  $\A (\Phi^{(k)} - \Phi) = \A \Phi^{(k)}-\Phi$. Thus 
 $\A \Phi^{(k)} \rightarrow \Phi$.  Hence $\A$ is continuous at $\Phi$.

Now consider $r \ge 1$ and suppose the assertion holds for $r-1$. Since $t_0 \neq 0$
the compact operator $H_\Phi$ is not zero and so  $H_\Phi^*  H_\Phi$ has
finite-dimensional eigenspace corresponding to $t_0^2$.  We now proceed as in the proof
of Theorem \ref{t1.1}: pick Schmidt vectors $v,~v^{(k)},~w,~w^{(k)}$, thematic functions
$V,~V^{(k)},~W^t,~W^{(k)t}$ and $L^\infty$ functions $Q,~Q^{(k)},~\Psi,~\Psi^{(k)}$
exactly as described above. Once again (\ref{recur1}) holds and the indices
corresponding to any nonzero superoptimal singular value in any thematic factorization
of $\Psi - \A\Psi$ are all 1.  Moreover, the superoptimal singular values of $\Psi$ are
$t_1,\dots,t_{m-1}$
, so that $\Psi$ has $r-1$ nonzero  superoptimal singular values.  By the inductive
hypothesis $\A\Psi^{(k)} \rightarrow \A\Psi$ in $X(M_{m-1,n-1})$.  The relations
(\ref{recur2}) now show that  $\A \Phi^{(k)} \rightarrow \Phi$ in $X(M_{m,n})$ as $k
\rightarrow \be$.  Thus  $\A$ is continuous at $\Phi$. \hfill $\bl$

\section{Necessary conditions for continuity}
\setcounter{equation}{0}
\setcounter{section}{2}
It is conceivable that the sufficient condition for continuity of $\A$ which we
established in Theorem \ref{t1.1} 
is also necessary for functions belonging to a space $X$
satisfying our axioms A1 to A5.  We can prove it for square
matrix functions whose superoptimal singular values are all nonzero.
\begin{lem}\thd
Let $\Phi \in X$ be of type $n \times n,$ and let $\e > 0$.  
Suppose that all $n$ superoptimal singular values of $\Phi$ are nonzero and that $\A$
is continuous at $\Phi$ with respect to the norm of $X$.  Then there
exists $\Psi \in X$ such that $\|\Phi -\Psi \|_X < \e$,
all $n$ superoptimal singular values of $\Psi$ are nonzero and all $n$ 
indices of $\Psi$ are equal to 1.
\end{lem}

\pf   Since $\A$ is continuous at $\Phi$ the same is true for the mapping $G
\mapsto \det(G-\A G)$, which maps $X(M_{n,n})$ to the space of constant functions
in $X$; it maps $G$ to the product of the superoptimal singular values of $G$.  
The latter mapping is nonzero at $\Phi$, by hypothesis, and
hence there exists $\e_1 > 0$ such that the product of the superoptimal
singular values of $G$ is nonzero whenever $\| \Phi - G \|_X < \e_1.$
It will therefore suffice to prove by induction on $n$ the following \\
\indent   Assertion:
{\em Let $\Phi \in X$ be of type $n \times n,$ and let $\e, \e_1 > 0$.  
Suppose that all $n$ superoptimal singular values of $G$ are nonzero
whenever $\| \Phi - G \|_X < \e_1.$  Then there
exists $\Psi \in X$ such that $\|\Phi -\Psi \|_X < \e$,
all $n$ superoptimal singular values of $\Psi$ are nonzero and all $n$ 
indices of $\Psi$ are equal to 1. }

To prove this we show first that there exists $\U \in X$ such that $\|H_\U \| > \|
H_\Phi\|$, $\|\U - \Phi\|_X$ is arbitrarily small
and $H_{\U-\Phi}$ has rank one.  Indeed, if $H_\Phi$ has maximising vector
$v$, $H_\Phi v = \bar{z} \bar{g}$ for some $g \in H^2$ and $\z \in \Bbb D$ is a point at
which $v$ is non-zero, then it suffices to take 
$$
\U(z) = \Phi(z) + (z-\z)^{-1} \eta \otimes v(\z)
$$
where $\eta \in \C^n$ is a non-zero vector of suitably small norm  satisfying
$\eta^t g(\z) > 0$.  We have
\bay
(H_\U v, \bar{z} \bar{g}) & = & 
((\Phi  + (z- \z)^{-1} \eta \otimes v(\z))v, \bar{z} \bar{g}) \nonumber \\
& = & (\Phi v, \bar{z} \bar{g})
+ (\pp_{-}  (z-\z)^{-1} \eta \otimes v(\z) v,\bar{z} \bar{g}) \nonumber \\
&=&  (H_\Phi v, \bar{z} \bar{g}) + \|v(\z)\|^2 ( (z-\z)^{-1}\eta, \bar{z} \bar{g})
= \|H_\Phi v \|^2 + \|v(\z)\|^2 \eta^t g(\z) \nonumber \\ 
&> &\|H_\Phi v \|^2 = \|H_\Phi \| ~\|v\| ~\|\bar{z} \bar{g} \|.\nonumber
\ey
Thus $\|H_\U \| > \| H_\Phi\|$. 
By choosing $\eta$ small we can ensure that $\U$ and $\Phi$ are close in any
norm, in particular the $X$ norm.  $\U$ thus has the
properties claimed. 

 Since $H_\U$ is a rank one perturbation of $H_\Phi$
$$
s_1(H_\U) \le s_0(H_\Phi) < s_0(H_\U),
$$
so that the maximising subspace of $H_\U$ is one-dimensional.
Let the superoptimal singular values of $\U$ be $t^\sharp_j,~j \geq 0$.
The index
of $t^\sharp_0 =  s_0(H_\U)$ in any thematic factorisation of $\U$ is 1;  for suppose otherwise.
Then we have
\begin{equation}
\label{induct}
\U - \A\U  = W^*  \left(\begin{array}{cc} t^\sharp_0 u & 0 \\ 0 & F \end{array} \right) V^*
\end{equation}
where $V,~W^t$ are thematic functions, $u$ is a badly approximable unimodular function 
and the Toeplitz 
operator $T_u$ has index less than --1.  Thus dim Ker $T_u > 1.$   It is easy to see that 
$\{Vf: f \in \Ker T_u \}$ is a space of maximising vectors of $H_\U$, 
and this contradicts the
simplicity of the singular value $s_0(H_\U)$.
  Thus the index of $t^\sharp_0$ is 1.  
Moreover, by [PY3, Theorem 1.1], 
$t^\sharp_1 \leq s_1(H_\U)$, so that $t^\sharp_1 < s_0(H_\U)$.
That is, in (\ref{induct}) $\|F\|_\be < t^\sharp_0$.

The case $n=1$ of Assertion is established by choice of $\Psi$ equal to $\U$.
Now consider $n > 1$ and suppose it true for $n-1$.   Pick $\U$ as above 
with  
$$
\|\U - \Phi\|_X <  \frac{1}{2} \min \{ \e, \e_1 \},
$$
and pick a thematic factorization (\ref{induct}) 
of $\U - \A\U$, so that $\|F\|_\be < t^\sharp_0$.
Since multiplication is continuous in the normed algebra $X(M_{n,n})$ there exists $K>1$
such that
$$
\|W^* G V^*\|_X \leq K \|G\|_X
$$
for all $G \in X(M_{n,n})$.  Let
$$
\d \df \min \{ \frac{\e}{2K}, \frac{\e_1}{2K}, t^\sharp_0 - \|F\|_\be \}.
$$
In the notation of (\ref{VW}) we have
$$
F = \b^* (\U - \A\U) \bar{\a} \in X(M_{n-1,n-1}).
$$
We claim that, for any $E \in  X(M_{n-1,n-1})$ such that $\|F - E \|_X < \d,$ the
superoptimal singular values of $E$ are all nonzero.  We have
\bay
\|E\|_\be &=& \|F\|_\be + \|E-F\|_\be \leq \|F\|_\be + \|E-F\|_X  \nonumber \\
& < & \|F\|_\be + t^\sharp_0 - \|F\|_\be = t^\sharp_0,  \nonumber
\ey
and hence
$$
\|  W^* \left( \begin{array}{cc} t^\sharp_0u & 0 \\ 0 & E \end{array} \right)V^* \|_\be = t^\sharp_0.
$$
Now let
\begin{equation}
\label{phiE}
\Phi_E = \A\U + W^* \left( \begin{array}{cc} t^\sharp_0u & 0 \\ 0 & E \end{array} \right)V^*.
\end{equation}
Then  $\|H_{\Phi_E}\| \leq \| \Phi_E - \A\U\|_\be = t^\sharp_0.$  Now $V, ~u$ have the form
$$
V= (v_{(i)} \quad \bar{\a}), \quad u=\bar{z} \bar{h}/ h
$$
where $v=v_{(i)} h$ is the inner-outer factorization of a maximising vector $v$
of $H_\Phi$ (see [PY1, Section 2, or PY2]).  This $v$ satisfies
$$
\|H_{\Phi_E}v \| = t^\sharp_0 \|v\|,
$$
and hence we have $\|H_{\Phi_E} \| = t^\sharp_0$.  Thus $\A\U$ is a 
best (though typically {\em
not} a superoptimal) analytic approximation to $\Phi_E$, and (\ref{phiE}) is a first stage
thematic factorization of $\Phi_E - \A\U$.  It follows from [PY1, Lemma 2.4] that
the superoptimal singular values of $E$ are those of  $\Phi_E$, all but the first.
However,
\bay
\label{error}
\|\Phi - \Phi_E \| & < & \|\Phi - \U\|_X + \|\U -\Phi_E\|_X
< \frac{1}{2} \min\{ \e, \e_1 \} + 
\| W^* \left( \begin{array}{cc} 0 & 0 \\ 0 & E - F\end{array} \right)V^*\|_X \nonumber \\
& \leq & \frac{1}{2} \min\{ \e, \e_1 \} + K\d <  \min\{ \e, \e_1 \}.
\ey
By hypothesis the superoptimal singular values of $\Phi_E$ are nonzero, and hence those
of $E$ are also. This establishes the claim.

 By the inductive hypothesis there exists $G \in X(M_{n-1,n-1})$ such
that $$\|F-G \|_X < \d,$$ 
all superoptimal singular values of $G$ are nonzero and all
$n-1$ indices of $G$ are 1.  Let
\begin{equation}
\label{psi}
\Psi \df \A\U + W^* \left( \begin{array}{cc} t^\sharp_0 u & 0 \\ 
0 & G\end{array} \right)V^*
\in X(M_{n,n}).
\end{equation}
In other words, $\Psi = \Phi_G$, and so by the above, the  superoptimal singular values
of $\Psi$ consist of $t^\sharp_0$ and those of $G$, 
hence are all nonzero.  By (\ref{error}),
$\|\Phi -\Psi \|_X < \e.$  Any thematic factorisation of $G-\A G$ induces one of
$\Psi - \A\Psi$ through the relation
$$
\Psi - \A\U -\b\A G\a^t =  W^* \left( \begin{array}{cc} t^\sharp_0 u & 0 \\ 0 & G -\A
G\end{array} \right)V^*,
$$
where we use the notation (\ref{VW}) for $V,~W$.  Since the indices of $t^\sharp_0 u$ and
$G-\A G$ are all 1, so are those of $\Psi - \A\Psi$.  The Assertion follows by
induction. $\hfill \bl$

\begin{thm} \thd
Let $X$ be a space of functions on $\T$ satisfying Axioms {\rm (A1)} to {\rm (A5)},
let $\Phi  \in X$ be of type $n \times n$ and suppose that the superoptimal singular
values of $\Phi$ are all nonzero.  If $\A$ is continuous at $\Phi$ then all indices in 
any thematic factorisation of $\Phi - \A \Phi$ are equal to 1.
\end{thm}

\pf Thematic functions have
constant determinant [PY1, Theorem 1.2].  
Hence $\det(\Phi - \A\Phi)$ is a function of nonzero
constant modulus on $\T$
whose winding number about $0$  is the sum of the indices in
any thematic factorisation of  $\Phi - \A\Phi$.  
Thus the winding number is $n$ if and only if
all the indices in any thematic factorisation are equal to 1.  By Lemma 2.1, 
$\Phi - \A\Phi$ is a limit in the norm of $X$ of a sequence of functions $\Psi$
such that $\Psi - \A\Psi$ has all indices defined and
equal to 1, hence such that $\det(\Psi -
\A\Psi)$ has winding number $n$.  It follows that  $\det(\Phi - \A\Phi)$ has winding
number $n.  \hfill \bl$

{\bf Remark.}  The proof shows a slightly stronger statement: if $\A$ is continuous at
$\Phi$ as a mapping from $X$ to $BMO$ (which is a weaker hypothesis than continuity
from $X$ to $X$) then the same conclusion holds.

As we mentioned in our discussion of sufficiency, continuity of $\A$ at functions 
which have some superoptimal singular value equal to zero is related to the
boundedness properties of scalar $\A$ on $X$. 
\begin{thm} \thd
Let $X$ be one of the Besov spaces $B^s_p,~s > 1/p$ or the Holder-Zygmund spaces 
$~\l_\a, ~\L_\a, ~\a > 0$.  Then $\A$ is
discontinuous at any matrix-valued function in $X$ which has a zero superoptimal
singular value.
\end{thm}
 
\pf It is shown in [P2] that $\A$ is unbounded on these spaces.
Let $\Phi \in X(M_{m,n})$.  We can suppose that $m \leq n$.  Let $t_r = 0,$ 
some $r \leq m$, but $t_j \neq 0$ for 
$j < r$. We suppose $r \geq 1$: the modifications for the case $r=0$ (i.e. $\Phi \in
H^\be$) are easy.  Consider a thematic factorisation
$$
\Phi - \A\Phi = W^*_0 \cdots W^*_{r-1} \left( \begin{array}{cccc} t_0u_0 & \cdot & 0& 0
\\ \cdot & \cdot & \cdot & \cdot \\  0 & \cdot &  t_{r-1} u_{r-1} & 0 \\
0 & \cdot & 0 & 0 \end{array} \right)  V_{r-1}^* \cdots V_0^*.
$$
By [P1], for $0 < \d < t_0$ we may pick a scalar function $\psi_\d \in X$ such that
$\|\psi_\d\|_X < \d$ and $\| \A \psi_\d \|_X \geq 1.$  Let
$$
\Phi_\d = \A \Phi +   W^*_0 \cdots W^*_{r-1}
\left( \begin{array}{ccccc} t_0u_0 & \cdot & 0& 0 & 0
\\ \cdot & \cdot & \cdot & \cdot  & \cdot\\  0 & \cdot &  t_{r-1} u_{r-1} & 0 & 0 \\
0 & \cdot & 0 & \psi_\d & 0 \\ 0 & \cdot & 0 & 0 & 0
\end{array} \right)  V_{r-1}^* \cdots V_0^*.
$$
Clearly $\| \Phi - \Phi_\d \|_X \to 0$ as $\d \to 0$.  If we solve the superoptimal
analytic approximation problem for $\Phi_\d$ by successive diagonalisation
then for the first $r$ stages it proceeds exactly as for $\Phi$ (a detailed proof of
this statement would be along the same lines as the proof of Lemma 2.1).  It follows that
$$
\Phi_\d - \A \Phi_\d =   W^*_0 \cdots W^*_{r-1}
\left( \begin{array}{ccccc} t_0u_0 & \cdot & 0& 0 & 0
\\ \cdot & \cdot & \cdot & \cdot  & \cdot\\  0 & \cdot &  t_{r-1} u_{r-1} & 0 & 0 \\
0 & \cdot & 0 & \psi_\d - \A \psi_\d& 0 \\ 0 & \cdot & 0 & 0 & 0
\end{array} \right)  V_{r-1}^* \cdots V_0^*.
$$
Thus
$$
\A\Phi - \A\Phi_\d =   W^*_0 \cdots W^*_{r-1}
\left( \begin{array}{ccccc} 0 & \cdot & 0& 0 & 0
\\ \cdot & \cdot & \cdot & \cdot  & \cdot\\  0 & \cdot &  0 & 0 & 0 \\
0 & \cdot & 0 & \A \psi_\d& 0 \\ 0 & \cdot & 0 & 0 & 0
\end{array} \right) V_{r-1}^* \cdots V_0^*.
$$
Since $\| \A \psi_\d\|_X \geq 1$, it cannot be true that $\A \Phi_\d \to \A\Phi$ in
$X$.  Thus $\A$ is discontinuous on $X$ at $\Phi. \hfill \bl$

\section{Continuity of superoptimal singular values}
The first superoptimal singular value $t_0$ of $\Phi \in H^\be + C$ is equal to
$\|H_\Phi\|$, hence is continuous with respect to the $L^\be$ norm.  Is the same true for
the other superoptimal singular values? Or at least with respect to one of the norms
$\| \cdot \|_X$ discussed above?  We will not venture a guess as to the answer to this
question, but we can at least prove continuity with respect to $\| \cdot \|_X$ 
under the same hypothesis as in Theorem 1.1.  For $\Phi \in H^\be+C$ we shall denote by
$t_j(\Phi)$ the $j$th superoptimal singular value of $\Phi$.
\begin{lem}\thd
Let $X \subset H^\be + C$ be a normed algebra of functions on $\T$ whose norm
majorises the $L^\be$ norm and which is invariant under $\A$.  If $\Phi \in X$ is a
point of continuity of $\A$ in $X$ then $\Phi$ is also a point of continuity of each of
the superoptimal singular values $t_j(\cdot)$ with respect to $\| \cdot \|$.
\end{lem}

\pf We recall that, for any matrix $A$ of type $m \times n$ and any integer $p$, 
$2 \leq p \leq m$, the $p$th exterior power $\wedge^p A$ is defined to be the
matrix of type
${m\choose p} \times {n\choose p}$ whose entries are the $p \times p$ minors of $A$.
Consider an $m \times n$ matrix function $G \in X, ~m \leq n$, and any integer $p$, 
$2 \leq p \leq m$. Define $(\wedge^p G)(z)$ to be  $\wedge^p(G(z))$.
Since the 
 entries of  $\wedge^p G$ are polynomials in those
of $G$ we have  $\wedge^p G \in X$ and the mapping $G \mapsto  \wedge^p G$ is continuous
with respect to the $X$ norms.  Thus, if $\A$ is continuous at $\Phi$, so is the mapping
$G \mapsto \| \wedge^p (G - \A G) \|_\be$.  It is immediate from consideration of
thematic factorisations that $\| \wedge^p (G - \A G) \|_\be$ equals the product of the
first $p$ superoptimal singular values of $G$.  Hence $t_0(\cdot),~t_0(\cdot)t_1(\cdot),
~t_0(\cdot)t_1(\cdot)t_2(\cdot),\cdots$ are all continuous at $\Phi$.  The result now
follows from the following simple observation which is valid for any topological space.
If $f_0 \geq f_1 \geq f_2 \geq \cdots \geq 0$ are real-valued functions such that
$f_0,~f_0f_1, ~f_0f_1f_2, \cdots$ are all
continuous at a point $x$ then each $f_j$ is continuous at $x$ (consider separately the
two cases $f_{j-1}(x) \neq 0$ and  $f_{j-1}(x) = 0$). $\hfill \bl$
\begin{thm}\thd
Let $X$ be a space of functions on $\Bbb T$ satisfying Axioms {\rm (A1)} to {\rm (A5)} 
and let  $\Phi\in X(M_{m,n})$, $m\leq n$.  Suppose that either $t_{m-1}\neq 0$ or $X$ is
the Besov space $B^1_1$.
If $\Phi-{\cal A}\Phi$ has a thematic factorisation with indices corresponding to
nonzero superoptimal singular values
all equal to 1
then  $t_j(\cdot)$ is continuous at $\Phi$  with respect to $\|\cdot \|_X$ for $0 \leq
j <m$.
\end{thm}
The proof is immediate from Theorems 1.1  and 1.10 and the foregoing Lemma.

\section*{References}

[AAK] {\sc V.M. Adamyan, D.Z. Arov and M.G. Krein}, Infinite Hankel block matrices and
some related continuation problems, {\em Izv. Akad. Nauk Armyan. SSR Ser. Mat}
{\bf 6} (1971), 87--112.

[F] {\sc B. Francis}, ``A Course in $H_\be$ Control Theory", Springer Verlag, Berlin, 1987.

[He] {\sc J. W. Helton}, The distance from a function to $H^\infty$ in the Poincar\'{e}
metric; electrical power transfer, {\em J. Functional Analysis,} {\bf 38} (1980), No.
2, 273--314.

[H] {\sc K. Hoffman}, ``Banach Spaces of Analytic Functions", Prentice Hall, Englewood
Cliffs, 1962.

[M] {\sc O. Merino}, Stability of qualitative properties and continuity of solutions to
problems of optimisation over spaces of analytic functions, Preprint 1989.

[Pa] {\sc M. Papadimitrakis}, Continuity of the operator of best approximation, {\em
Bull. London Math. Soc.} {\bf 25} (1993)44-48.

[Par] {\sc J. R. Partington}, Robust identification and and interpolation in $H^\be$,
{\em Int. J. Control} {\bf 54}(1991) 1281-1290.

[P1] {\sc V.V.Peller}, Hankel operators and continuity properties of
best approximation operators, {\em Algebra i Analiz}, {\bf 2}:1, (1990),
163-189. English Transl. in {\em Leningrad Math. J.}, {\bf 2} (1991), 139-160.

[P2] {\sc V.V. Peller}, Boundedness properties of the operators of best approximation by
analytic and meromorphic functions, {\em Ark. Mat} {\bf 30} (1992) 331-343.

[P3] {\sc V.V. Peller}, Approximation by analytic operator-valued functions

[PK] {\sc V.V. Peller and S.V. Khruschev}, Hankel operators, best approximation and
stationary Gaussian processes, {\em Russian Math. Surveys,} {\bf 37} (1982)
53--124.

[PY1] {\sc V.V. Peller and N.J. Young}, Superoptimal analytic approximations of matrix
functions, {\em J. Functional Analysis} {\bf 120}(1994) 300-343.

[PY2] {\sc V.V. Peller and N.J. Young}, Construction of superoptimal
 approximants, {\it Mathematics of Signals, Systems and Control}, to appear.

[PY3] {\sc V.V. Peller and N.J. Young}, Superoptimal singular values and indices of
matrix functions, {\it Integral Equations and Operator Theory} {\bf 20} (1994) 350-363.

[T1] {\sc S.R. Treil}, The Adamyan-Arov-Krein theorem: a vector version, {\em
Zap. Nauchn. Semin. LOMI} {\bf 141} (1985), 56--71 (Russian).

[T2] {\sc S.R. Treil}, On superoptimal approximation by analytic and meromorphic
matrix-valued functions, {\em J. Functional Analysis,} {\bf 131} (1995) 386-414.

[Yo] {\sc K. Yosida}, ``Functional Analysis" (Sixth edition), Springer Verlag, Berlin
1980.

 \end{document}